\documentclass[a4paper,10pt,reqno]{amsart}
\usepackage{enumitem}
\usepackage{amsmath}
\usepackage{amsthm}
\usepackage{amssymb}
\usepackage[capitalize]{cleveref}

\newcounter{alph}

\newtheorem{theo}[alph]{Theorem}

\newtheorem{coro}[alph]{Corollary}
\numberwithin{equation}{section}
\newtheorem{cor}[equation]{Corollary}
\newtheorem{lem}[equation]{Lemma}

\newtheorem{prop}[equation]{Proposition}
\newtheorem{thm}[equation]{Theorem}
\theoremstyle{definition}

\newtheorem{exa}[equation]{Example}
\newtheorem{rem}[equation]{Remark}

\DeclareMathOperator{\im}{im}

\DeclareMathOperator{\diam}{diam}

\DeclareMathOperator{\ess}{ess}

\DeclareMathOperator{\Ray}{Ray}

\DeclareMathOperator{\Vol}{Vol}

\def\R{\mathbb R}
\def\Z{\mathbb Z}

\def\ve{\varepsilon}
\def\vf{\varphi}
\def\la{\langle}
\def\ra{\rangle}

\begin{document}

\title[Spectral instability]{Spectral instability of coverings}
\author{Werner Ballmann}
\address{WB: Max Planck Institute for Mathematics, Vivatsgasse 7, D53111 Bonn}
\email{hwbllmnn\@@mpim-bonn.mpg.de}
\author{Sugata Mondal}
\address{SM: Max Planck Institute for Mathematics, Vivatsgasse 7, D53111 Bonn,
and Department of Mathematics and Statistics, University of Reading, UK}
\email{s.mondal@reading.ac.uk}

\date{\today}


\thanks{\emph{Acknowledgments.}
We are grateful to the Max Planck Institute for Mathematics
and the Hausdorff Center for Mathematics in Bonn for their support and hospitality.
The second named author would like to thank C.S.\,Rajan and C.\,Drutu for helpful discussions.}

\begin{abstract}
We study the behaviour of eigenvalues, below the bottom of the essential spectrum, of the Laplacian under finite Riemannian coverings of complete and connected Riemannian manifolds.
We define spectral stability and instability of such coverings.
Among others, we provide necessary conditions for stability or, equivalently, sufficient conditions for instability.
\end{abstract}

\maketitle


\setcounter{tocdepth}{2}

\section{Introduction}

Recently Magee et al.\;\cite{MageeNaud20,HildeMagee21} have initiated a study of the spectrum of the Laplacian of a random Riemannian cover of a fixed hyperbolic (i.e., curvature $=-1$) surface. 
Broadly speaking, the main results obtained by them say that asymptotically almost surely (with respect to the uniform measure on the space of $n$-sheeted covers) the spectrum of a covering surface does not acquire a new eigenvalue below a specific threshold $<1/4$.
What is even more interesting is that this threshold is independent of the bottom surface,
but only depends on its type; see below.

A classical result of Randol \cite{Ra74} is quite opposite to the results of Magee et al.
Namely, for any hyperbolic metric on a \emph{closed} (i.e., \emph{compact and connected with empty boundary}) surface $S$ of genus $g\ge2$,
any natural number $\ell$ and any $\ve>0$,
there is a finite Riemannian covering $p\colon S'\to S$ such that $S'$ has at least $\ell$ eigenvalues in $[0,\ve)$.
(See \cite[Theorem 4.1]{BMM18} for an elementary proof.)

One motivation for our studies in this paper is that, although the above mentioned results of Magee et al.\;show 
that asymptotically almost all $n$-sheeted covers of the surface in question are {\it spectrally stable} in a specific range,
they do not provide any necessary or sufficient condition for this to happen.
Among others, we provide, in this paper, some necessary conditions.

To set the stage, let $M$ be a complete and connected Riemannian manifold of dimension $m$.
Denote by  $\tilde M$ the universal covering space of $M$, endowed with the lifted Riemannian metric,
and let $\Gamma$ be the fundamental group of $M$,
viewed as the group of covering transformations on $\tilde M$.

Denote by $\lambda_0(M)\le\lambda_{\ess}(M)$ the bottom of the spectrum and the essential spectrum (of the Laplacian $\Delta$) of $M$, respectively.
Recall that $\lambda_0(M)$ need not vanish in general,
that $\lambda_0(M)=0$ and $\lambda_{\ess}(M)=\infty$ in the case where $M$ is closed
(that is, compact and connected without boundary)
and that the spectrum of $M$ below $\lambda_{\ess}(M)$ consists of a locally finite set of eigenvalues of finite multiplicity.
We assume throughout that
\begin{align}\label{zess}
	\lambda_0(M) < \lambda_{\ess}(M)
\end{align}
and enumerate the eigenvalues of $M$ in $[0,\lambda_{\ess}(M))$ according to their size and multiplicity as
\begin{align}\label{enum}
    0 \le \lambda_0(M) < \lambda_1(M) \le \lambda_2(M) \le \dots,
\end{align}
where we recall that the eigenvalue $\lambda_0(M)$ has multiplicity one since its eigenfunctions do not change sign.
For $\lambda\ge0$, we denote by $N_M(\lambda)$ and $N_M(\lambda-)$ the number of $\lambda_k(M)$ in $[0,\lambda]$ and $[0,\lambda)$, respectively.
More generally, for any subset $I\subseteq[0,\lambda_{\ess}(M))$,
we denote by $N_M(I)$ the number of $\lambda_k(M)$ in $I$.

Let $p\colon M'\to M$ be a finite Riemannian covering of complete and connected Riemannian manifolds.
Then
\begin{align}\label{bott}
    \lambda_0(M)=\lambda_0(M')
    \quad\text{and}\quad
    \lambda_{\ess}(M)=\lambda_{\ess}(M'),
\end{align}
see \eqref{bottomeq}.
Since the lifts $p^*\vf=\vf\circ p$ of eigenfunctions $\vf$ of $M$ are eigenfunctions of $M'$,
we always have
\begin{align}\label{liftk}
	\lambda_k(M') \le \lambda_k(M).
\end{align}
Likewise, for any subset $I\subseteq[0,\lambda_{\ess}(M))$,
\begin{align}\label{lifti}
     N_M(I) \le N_{M'}(I).
\end{align}
We say that $p$ is \emph{$I$-stable} if 
\begin{align}\label{stabi}
     N_M(I) = N_{M'}(I).
\end{align}
With respect to this terminology, the results of Magee et al.\;say that, \\
1) for any orientable, convex cocompact, non-compact hyperbolic surface $S$ with Hausdorff dimension of its limit set $\delta>1/2$ and any $\sigma\in(3\delta/4,\delta)$,
any finite Riemannian cover $p\colon S'\to S$ is asymptotically almost surely $[\delta(1-\delta),\sigma(1-\sigma)]$-stable as $|p|\to\infty$;
see \cite{MageeNaud20}.
Note that here $\lambda_0(S)=\delta(1-\delta)<1/4=\lambda_{\ess}(S)$. \\
2) for any orientable and compact hyperbolic surface $S$ and any $\ve>0$,
any finite Riemannian cover $p\colon S'\to S$ is  asymptotically almost surely $[0,3/16-\ve]$-stable as $|p|\to\infty$;
see \cite{MageeNaudPuder22}. \\
3) for any orientable and non-compact hyperbolic surface $S$ of finite area and any $\ve>0$,
any finite Riemannian cover $p\colon S'\to S$ is asymptotically almost surely $[0,1/4-\ve]$-stable as $|p|\to\infty$;
see \cite{HildeMagee21}.
Here $\lambda_0(S)=0<1/4=\lambda_{\ess}(S)$.

Clearly, $I$-stability means that lifting yields an isomorphism between eigenspaces of $M$ and $M'$ for all eigenvalues of $M$ and $M'$ in $I$.
In particular, \\
1) if $J\subseteq[0,\infty)$ is a further subset and $I\subseteq J$, then
\begin{align}\label{IJ}
	\text{$I$-instability of $p$ implies $J$-instability of $p$.}
\end{align}
2) if $q\colon S''\to S'$ is a further finite Riemannian covering of complete and connected Riemannian manifolds, then
\begin{align}\label{pq}
	\text{$I$-instability of $p$ or $q$ implies $I$-instability of $p\circ q$.}
\end{align}
For $k>0$, we say that $p$ is $\lambda_k$-\emph{stable} if
\begin{align}\label{stabk}
     N_M(\lambda) = N_{M'}(\lambda),
\end{align}
where $\lambda=\lambda_k(M)<\lambda_{\ess}(M)$.
By definition, $\lambda_k$-instability $N_M(\lambda) < N_{M'}(\lambda)$  can occur in two ways:
Either $p$ is \emph{strictly $\lambda_k$-unstable}, that is, $N_M(\lambda-)<N_{M'}(\lambda-)$,
or else $p$ is \emph{weakly $\lambda_k$-unstable}, that is, $N_M(\lambda-)=N_{M'}(\lambda-)$,
but the multiplicity of $\lambda$ as an eigenvalue increases, $\mu(\lambda,M)<\mu(\lambda,M')$.
By \eqref{IJ},
\begin{align}\label{kell}
	\text{$\lambda_k$-instability implies $\lambda_\ell$-instability, for any $1\le k\le\ell$,}
\end{align}
where $\lambda_\ell(M)<\lambda_{\ess}(M)$.
In particular, $\lambda_1$-instability implies $\lambda_k$-instability for all $k\ge1$.
For that reason, our main focus is on $\lambda_1$-stability and instability.

For an eigenfunction $\vf$ on a Riemannian manifold,
the set $\mathcal{Z}(\vf)=\{\vf=0\} $ is called the \emph{nodal set of $\vf$}
and the connected components of $\{\vf\ne0\}$ are called \emph{nodal domains of $\vf$}.
(In general, we use the term \emph{domain} to indicate \emph{connected open sets}.)
One of our main arguments uses connectedness of preimages in $M'$ of nodal sets in $M$.

\begin{theo}\label{thma}
If $p\colon M'\to M$ is a $\lambda_k$-stable finite Riemannian covering of complete and connected Riemannian manifolds, where $\lambda_0(M)<\lambda_k(M)<\lambda_{\ess}(M)$,
then the preimage $p^{-1}(U)$ of any nodal domain $U$ of any $\lambda$-eigenfunction $\vf$ on $M$ is connected, for any $\lambda_0(M)\le\lambda\le\lambda_k(M)$.
In fact, if $U$ is any nodal domain of any $\lambda$-eigenfunction $\vf$ on $M$ for any such $\lambda$
and $j\ge1$ denotes the number of connected components of $p^{-1}(U)$, then
\begin{align*}
	N_{M'}(\lambda-) \ge N_{M}(\lambda-) + j - 1.
\end{align*}
\end{theo}

\cref{thma} is a special case of \cref{numberg}.
An easy application of the main results of \cite{Ko06} and \cite{EP15} and \cref{thma} yields the following

\begin{coro}\label{pball}
A closed manifold $M$ carries a Riemannian metric $g$,
such that $p$ is strictly $\lambda_1$-unstable with respect to $g$,
for any non-trivial finite covering $p\colon M'\to M$ of closed and connected manifolds,
where $M'$ is endowed with the lifted Riemannian metric $g'$.
In fact, for an appropriate choice of $g$,
\begin{align*}
	N_{M'}(\lambda-) \ge |p| \quad(= N_{M}(\lambda-) + |p| - 1),
\end{align*}
where $\lambda=\lambda_1(M,g)$ and $|p|$ denotes the degree (number of sheets) of $p$.
\end{coro}

\begin{proof}
By \cite[Main Theorem]{Ko06} and \cite[Theorem 1.1]{EP15}, $M$ carries a Riemannian metric $g$,
which has a topological ball $U$ as a nodal domain.
(Note that the proof in \cite{Ko06} also works in the non-orientable case.)
Since balls are simply connected, $p^{-1}(U)$ has $|p|$ disjoint lifts.
Now the claim follows from \cref{thma}.
\end{proof}

\begin{coro}\label{examon}
The orientable closed surface $S$ of genus two carries a hyperbolic metric
such that any non-trivial Riemannian covering $p\colon S'\to S$,
that is not generated by one element, is strictly $\lambda_1$-unstable.
\end{coro}

Here we say that a covering $p\colon M'\to M$ of connected manifolds is \emph{generated by $k$ elements}
if, for one--or any--point $x\in M$, there are $k$ loops in $M$ at $x$
such that any two points in $p^{-1}(x)$ can be connected by lifts to $M'$ of concatenations of these loops and their inverses; see also \cref{subgen}.
This property is independent of the choice of $x$.

\begin{proof}[Proof of \cref{examon}]
By \cite{Mondal15}, there exists a hyperbolic metric on $S$ which has a $\lambda_1(S)$-eigenfunction $\vf$
such that one of its nodal domains $U$ is either a disc or an annulus.
In the first case, the preimage of $U$ is disconnected with $|p|$ components
and the assertion is a consequence of \cref{thma}.
In the second csae, if the preimage of $U$ is connected and $x\in U$,
then any two points of $p^{-1}(x)$ can be connected by a lift of an iterate of any loop in $U$ at $x$ which generates $\pi_1(U,x)$;
see also \cref{preu}.
This shows the assertion in the second case.
\end{proof}

\begin{theo}\label{theoi}
Suppose that $M$ is complete and connected with $\lambda_1(M)<\lambda_{\ess}(M)$ and carries a $\lambda_1(M)$-eigenfunction $\vf$
such that its nodal set is not connected.
\begin{enumerate}
\item Then there is a two-sheeted Riemannian covering of $M$ which is strictly $\lambda_1$-unstable.
\item If $M$ is orientable, then $M$ carries an $n$-sheeted cyclic Riemannian covering which is strictly $\lambda_1$-unstable, for any $n\ge2$.
\end{enumerate} 
\end{theo}

\cref{theoi} is a special case of \cref{theoi2}.

Let $S$ be a complete and connected Riemannain surface.
(To avoid misunderstandings: a Riemannian surface is a surface with a Riemannian metric.)
Recall that $S$ is said to be of \emph{finite type} if it is diffeomorphic to the interior of a compact surface with (possibly empty) boundary.

For a domain $U$ in $S$ and a point $x\in U$,
we identify $\Gamma$ with $\pi_1(S,x)$ and denote the \emph{image of $\pi_1(U,x)$ in $\Gamma$} by $\Gamma_U$.
Corresponding assertions will be independent of the choice of $x$.

\begin{theo}\label{thmb}
Assume that $S$ is of finite type with $\chi(S)<0$, and let $\vf$ be a $\lambda$-eigenfunction,
where $\lambda_0(S)<\lambda<\lambda_{\ess}(S)$.
Then $\vf$ has $\nu\ge2$ nodal domains and at least one, $U$, such that $\chi(S)/\nu\le\chi(U)\le1$.
For any such $U$, $\Gamma$ admits a surjective homomorphism $I$ to $\Z^{\mu}_2$,
respectively $\Z^{\mu}$ if $S$ is orientable,
where $\mu\ge-(\nu-1)\chi(S)/\nu$, such that $\Gamma_U\subseteq\ker I$.
In particular, if $\Gamma'\subseteq\Gamma$ is a finite index subgroup containing $\Gamma_U$,
then the corresponding Riemannian covering $p\colon S'\to S$ is strictly $\lambda$-unstable.
More precisely,
\begin{align*}
	N_{S'}(\lambda-) \ge N_S(\lambda-) + |p| -1,
\end{align*}
where $|p|=|\Gamma'\backslash\Gamma|$.
\end{theo}

\cref{thmb} is a special case of \cref{thmb2}.
It applies, for example, to non-compact hyperbolic surfaces $S$ of finite area
with $0<\lambda<1/4$ since, for them, $\lambda_{\ess}(S)=1/4$.
The number $\mu$ is determined in the proof of \cref{thmb2}.

\begin{rem}[Weyl's law]\label{weyl}
Let $p\colon M' \to M$ be a non-trivial finite Riemannian covering of closed Riemannian manifolds.
Then, by Weyl's  law,
\begin{align*}
	\lim_{\lambda\to\infty}\frac{N_M(\lambda)}{\lambda^{m/2}}=C_m\Vol M
	\quad\text{and}\quad
	\lim_{\lambda\to\infty}\frac{N_{M'}(\lambda)}{\lambda^{m/2}}=C_m\Vol M',
\end{align*}
where $C_m$ equals the volume of the ball of radius $1/2\pi$ in $\R^m$.
Since $\Vol M'=|p|\Vol M$ and $|p|\ge2$,
we get that $N_{M'}(\lambda)>N_M(\lambda)$ for all sufficiently large $\lambda$.
Therefore stability of $p$ can only hold in a bounded range of $\lambda$.
\end{rem}

A sequence of Riemannian coverings of Riemannian manifolds,
\begin{align*}
    \cdots \to M_k \to M_{k-1} \to \dots \to M_1 \to M_0 = M,
\end{align*}
is called a \emph{tower of Riemannan coverings}.
If the $M_k$ are connected, then the universal covering $\tilde M$ of $M$ sits above the tower.
We endow it with the lifted Riemannian metric.
The next result asserts that we should not expect stability beyond the bottom of the spectrum of $\tilde M$.

\begin{theo}\label{tower}
For a tower of finite Riemannian coverings of complete and connected Riemannian manifolds,
suppose that the degree of the coverings $M_k\to M$ tends to infinity and that $\lambda_0(\tilde M)<\lambda_{\ess}(M)$.
Then, for any $\lambda_0(\tilde M)<\lambda<\lambda_{\ess}(M)$ and $l\ge1$,
$M_k$ has at least $l$ eigenvalues below $\lambda$, for all sufficiently large $k$; in short,
\begin{align*}
	\limsup_{k\to\infty}\lambda_l(M_k) \le \lambda_0(\tilde M).
\end{align*}
\end{theo}

\cref{tower} is a special case of \cref{roof}.
It is motivated by \cite[Proposition 6]{Sunada85} and \cite[Theorem 2]{Brooks86}.
Note that we do not need residual finiteness of $\Gamma$ respectively that $\cap\Gamma_k=\{1\}$,
since our proof relies on the push down construction from \cite{BMM18}.

The assumption $\lambda_0(\tilde{M})<\lambda_{\ess}(M)$ is satisfied if $M$ is compact since the essential spectrum of closed Riemannian manifolds is empty.
On the other hand, if $\lambda_0(\tilde{M})=\lambda_{\ess}(M)$, the argument in the proof of \cref{tower} still applies,
but the assertion might be meaningless for any $\lambda>\lambda_{\ess}(M)$.
At least the spectral projection corresponding to $[0,\lambda)$ would then have infinite rank.

\section{Setup and preliminaries}
Let $M$ be a complete and connected Riemannian manifold of dimension $m$.
Let $\Delta$ denote the Laplace-Beltrami operator,
acting on the space of smooth functions $C^\infty(M)$ on $M$. 
Recall that $\Delta$ is {\em essentially self-adjoint} on $C^\infty(M)\subseteq L^2(M)$.
Its closure will also be denoted by $\Delta$.
It has domain $H^1(M)$, and its spectrum, depending on the context denoted by
\begin{align*}
    \sigma(M,\Delta) = \sigma(\Delta) = \sigma(M),
\end{align*}
can be decomposed into two sets,
\begin{align*}
    \sigma(M) = \sigma_d(M) \sqcup \sigma_{\ess}(M),
\end{align*}
the {\em discrete spectrum} and the {\em essential spectrum}.
Recall that $\sigma_d(M)$ consists of isolated eigenvalues of $\Delta$ of finite multiplicity
and that $\sigma_{\ess}(M)$ consists of those $\lambda\in\R$ for which $\Delta-\lambda$ is not a Fredholm operator.
By elliptic regularity theory, $\sigma(M)=\sigma_d(M)$ if $M$ is compact.
By the above characterization of the discrete spectrum, $\sigma(M)=\sigma_{\ess}(M)$ if $M$ is homogeneous and non-compact.

Denote by $\lambda_0(M)\le\lambda_{\ess}(M)$ the bottom of $\sigma(M)$ and $\sigma_{\ess}(M)$, respectively.
If $M$ is compact, then $\lambda_0(M)=0<\lambda_{\ess}(M)=\infty$.
Furthermore, $0$ is an eigenvalue of $\Delta$ of multiplicity one with constant functions as eigenfunctions.
In general, $\lambda_0(M)$ may be positive and may belong to $\sigma_{d}(M)$
or we may have $\lambda_0(M)=\lambda_{\ess}(M)$.
We shall be interested in the case where
\begin{align}\label{standardass}
    \lambda_0(M) < \lambda_{\ess}(M).
\end{align}
Then $\lambda_0(M)$ is an eigenvalue of $\Delta$ of multiplicity one with eigenfunctions which do not change sign.
Moreover, by the above, we have
\begin{align}\label{discretev}
    \sigma(M) \cap [0,\lambda_{\ess}(M)) \subseteq \sigma_{d}(M).
\end{align}
We enumerate the eigenvalues of $\Delta$ in $[0,\lambda_{\ess}(M))$ by size,
\begin{align}\label{enumev}
    \lambda_0(M) < \lambda_1(M) \le \lambda_2(M) \le \dots < \lambda_{\ess}(M),
\end{align}
where repetitions account for multiplicities.
In general, the number of eigenvalues of $\Delta$ below $\lambda_{\ess}(M)$ might be infinite, even if $\lambda_{\ess}(M)<\infty$.
On the other hand, we have the \emph{variational characterization} of $\lambda_k(M)<\lambda_{\ess}(M)$ by
\begin{align}
    \lambda_k(M) = \inf_F\max_{0\ne\vf\in F}\Ray\vf,
\end{align}
where $F$ runs over all subspaces of $H^1(M)$ of dimension $k+1$ and $\Ray\vf$ denotes the \emph{Rayleigh quotient} of $\vf$.
The infimum is achieved by the linear span of the $\lambda_j(M)$-eigenfunctions, where $0\le j\le k$.

\subsection{Spectrum under finite Riemannian coverings}\label{subcov}
Unless otherwise specified, Riemannian manifolds are assumed to be complete and connected.
Similarly, $p\colon M'\to M$ will denote a finite Riemannian covering of Riemannian manifolds with $|p|$ sheets and group $\Gamma$ of covering transformations.
Recall that $\Gamma$ is transitive on the fibers of $p$ if and only if $p$ is normal.
Since $p$ is finite, we have
\begin{align}\label{bottomeq}
    \lambda_0(M')=\lambda_0(M) \quad\text{and}\quad \lambda_{\ess}(M')=\lambda_{\ess}(M).
\end{align}
To show the first equality, recall that $\lambda_0$ is the supremum of the \emph{positive spectrum},
that pull back and averaging are inverses to each other on the respective spaces of positive functions on $M$ and $M'$,
and that both, pull back and averaging, are compatible with the Laplacian; see \cite[Theorem 2.1]{Su87}.
As for the second equality, recall the well known fact that $\lambda_{\ess}$ is the $\limsup$ of $\lambda_0$ on the family of neighborhoods of infinity
of the corresponding manifold and that the first equality does not require completeness; see e.g.\;\cite[Proposition 4.8]{BP21}.

As indicated in \eqref{standardass}, we assume throughout that
\begin{align}\label{zeroless}
    \lambda_0(M) < \lambda_{\ess}(M).
\end{align}
Then $\lambda_0(M)$ is an eigenvalue of $\Delta$ on $M'$ and $M$ of multiplicity one with unique positive eigenfunctions $\vf_0'$ and $\vf_0$ of respective $L^2$-norms equal to one.

For any functions $\vf'$ on $M'$ and $\vf$ on $M$,
let $p_*\vf'$ be the function on $M$ such that $(p_*\vf')(x)$ is the average of the $\vf'(y)$, $y\in p^{-1}(x)$,
and $p^*\vf=\vf\circ p$ be the pull-back of $\vf$ to $M'$.
Say that $\vf'$ is $p$-invariant if $\vf'$ is constant along the fibers of $p$.
Obviously, this holds if and only if there is a function $\vf$ on $M$ such that $\vf'=p^*\vf$ or, equivalently, if and only if $\vf'=p^*p_*\vf'$.
Clearly, $p_*$ and $p^*$ preserve all the standard regularity and integrability conditions.
For $\vf'\in L^2(M')$ and $\vf\in L^2(M)$, we have
\begin{align}\label{padjoint}
    \la\vf',p^*\vf\ra_{M'} = |p|\la p_*\vf',\vf\ra_M,
\end{align}
where the indices $M'$ and $M$ indicate the scalar products in $L^2(M')$ and $L^2(M)$, respectively.
Furthermore,
\begin{align}
    L^2(M') = \im p^* \oplus \ker p_*.
\end{align}
For any $\lambda\ge0$, let $E_\lambda'$ and $E_\lambda$ be the $\lambda$-eigenspaces of $\Delta$ on $M'$ and $M$, respectively.
Since lifts of $\lambda$-eigenfunctions on $M$ are $\lambda$-eigenfunctions on $M'$,
$p^*E_\lambda$ is equal to the space of $p$-invariant functions in $E_\lambda'$.

\begin{prop}\label{pstar}
For any $\lambda\ge0$,
\begin{enumerate}
    \item $E_\lambda'=p^*E_\lambda\oplus\ker p_*$ is an $L^2$-orthogonal splitting;
    \item $\sqrt{|p|}p_*\colon p^*E_\lambda\to E_\lambda$ is an orthogonal isomorphism with inverse $p^*/\sqrt{|p|}$.
\end{enumerate}
\end{prop}

Note here that $\lambda$-eigenspaces of $\Delta$ are closed subspaces of $L^2(M)$ and $L^2(M')$,
even if $\lambda$ belongs to the essential spectrum of $M$ respectively $M'$.

\section{Connectedness under coverings}\label{secco}
Fix a point $x\in M$.
For $x'\in p^{-1}(x)$ and a loop $c\colon[0,1]\to M$ at $x$,
let $c_{x'}$ be the lift of $c$ to $M'$ starting at $x'$.
Then $x'[c]=c_{x'}(1)$ defines a right action of $\Gamma=\pi_1(M,x)$ on the fiber $p^{-1}(x)$ of $p$ over $x$,
where $[c]\in\Gamma$ denotes the homotopy class of $c$.

\begin{lem}\label{preu}
Let $U\subseteq M$ be a connected open subset containing $x$.
Then the connected components of $p^{-1}(U)$ are in canonical one-to-one correspondence
with the orbits of $\Gamma_U$ on $p^{-1}(x)$,
where $\Gamma_U$ denotes the image of $\pi_1(U,x)$ in $\Gamma$. 
\end{lem}

\begin{proof}
Let $c\colon[0,1]\to U$ be a loop at $x$, $x'\in p^{-1}(x)$, and $c_{x'}$ be the lift of $c$ starting at $x'$.
Then the image of $c_{x'}$ is contained in $p^{-1}(U)$,
and hence $x'[c]$ belongs to the component of $p^{-1}(U)$ containing $x'$.
Conversely, if $x''\in p^{-1}(x)$ belongs to the same component of $p^{-1}(U)$ as $x'$,
then there is a path $c'\colon[0,1]\to p^{-1}(U)$ from $x'$ to $x''$.
Then $c'=c_{x'}$ and, therefore, $x''=c_{x'}(1)$, where $c=p\circ c'$ is a loop at $x$. 
\end{proof}

Fixing a point $x'\in p^{-1}(x)$, we get a canonical identification $p^{-1}(x)=\Gamma'\backslash\Gamma$,
where $\Gamma'$ denotes the image of $\pi_1(M',x')$ in $\Gamma$.
With respect to this identfication, the right action of $\Gamma$ on $p^{-1}(x)$
corresponds to the right action of $\Gamma$ on $\Gamma'\backslash\Gamma$.

\begin{cor}\label{preu2}
After the choice of a point $x'\in p^{-1}(x)$,
the connected components of $p^{-1}(U)$ are in canonical one-to-one correspondence
with the elements of $\Gamma'\backslash\Gamma/\Gamma_U$,
the space of orbits of the right action of $\Gamma_U$ on $\Gamma'\backslash\Gamma$.
\end{cor}

When it comes to the existence of $\lambda_k$-unstable coverings,
the case $\Gamma_U\subseteq\Gamma'$ is of interest.
Our next result corresponds to the lifting property of covering projections.

\begin{lem}\label{gugp}
If $\Gamma_U\subseteq\Gamma'$, then the right action of $\Gamma_U$ on $\Gamma'\backslash\Gamma$ fixes $\Gamma'e$.
Hence the action has more than one orbit unless $\Gamma'=\Gamma$.
If the normalizer $N_\Gamma(\Gamma_U)$ of $\Gamma_U$ is contained in $\Gamma'$, then the right action of $\Gamma_U$ on $\Gamma'\backslash\Gamma$ is trivial.
In this case, the action has $|p|=|\Gamma'\backslash\Gamma|$ orbits.
\end{lem}

\begin{proof}
If $g\in\Gamma_U$, then $\Gamma'eg=\Gamma'g=\Gamma'$ since $g\in\Gamma'$.
Under the second assumption, if $g\in\Gamma_U$ and $h\in\Gamma$,
then $\Gamma'hg=\Gamma'g'h=\Gamma'h$ since $g'\in N_\Gamma(\Gamma_U)\subseteq\Gamma'$.
\end{proof}

\begin{exa}[Abelian coverings]\label{abco}
For a domain $U\subseteq M$ and a point $x\in U$, consider the Hurewicz homomorphism
\begin{align*}
	H_x\colon \pi_1(M,x) \to \pi_1(M,x)/[\pi_1(M,x),\pi_1(M,x)] = H_1(M) 
\end{align*}
and the projection
\begin{align*}
	H_1(M) \to H_1(M)/i_*(H_1(U))=:A,
\end{align*}
where $H_1$ indicates first homology groups with coefficients in $\Z$
and $i\colon U\to M$ denotes the inclusion.
Under their composition, the preimage of $0\in A$ in $\pi_1(M, x)$ equals $\Gamma_U$.
Hence the preimage $\Gamma'$ of any finite index subgroup $A'$ of $A$ is a normal subgroup of $\Gamma$ containing $\Gamma_U$
such that $\Gamma'\backslash\Gamma\cong A'\backslash A$ is a finite Abelian group.
A question, among others adressed in \cref{theoi2} and \cref{thmb2}, is whether $A$ is trivial.
\end{exa}

\subsection{Minimal number of generators}\label{subgen}
Say that $\Gamma'\backslash\Gamma$ is \emph{generated by $k$ elements}
if there is a subset $G$ of $\Gamma$ with $|G|=k$ such that $\Gamma'\cup G$ generates $\Gamma$;
then the elements of $G$ are also called \emph{generators of $\Gamma'\backslash\Gamma$}.
In the case where $\Gamma'$ is a normal subgroup of $\Gamma$,
this terminology coincides with the usual one for the group $\Gamma'\backslash\Gamma$.
The \emph{minimal number of generators} of $\Gamma'\backslash\Gamma$ is denoted by $\mu(\Gamma'\backslash\Gamma)$.

\begin{lem}\label{preu4}
If $\mu(\Gamma'\backslash\Gamma)=k$ and $\Delta\subseteq\Gamma$ is a subgroup generated by $\ell$ elements,
then the right action of $\Delta$ on $\Gamma'\backslash\Gamma$ has at least $k-\ell+1$ orbits.
\end{lem}

\begin{proof}
The claim is true for $\ell\ge k$.
Suppose now, by induction, that it is true for $\ell+1\le k$.
Then the right action of the subgroup $\Delta_g$ of $\Gamma$ generated by $\Delta$ and any additional element $g\in\Gamma$
has at least $k-\ell$ orbits in $\Gamma'\backslash\Gamma$.
If all of these would be orbits of $\Delta$ already, for any choice of $g\in\Gamma$, then the $\Delta$-orbits would be invariant under $\Gamma$.
However, that cannot be because $\Gamma$ has only one orbit.
Hence there is a choice of a $g\in\Gamma$ which decreases the number of orbits by one. 
\end{proof}

\begin{rem}\label{mga}
In general, the calculation of the minimal number of generators of quotients $\Gamma'\backslash\Gamma$ is a difficult problem.
However, in the case where $\Gamma'$ is a normal subgroup of $\Gamma$ such that $A=\Gamma'\backslash\Gamma$ is a finite Abelian group,
\begin{align*}
	A = \mathbb{Z}/k_1\times\dots\times\mathbb{Z}/k_n,
\end{align*}
then $\mu(A)$ equals the maximal number of $k_i$ which share a common divisor.
\end{rem}

\subsection{Asymptotic estimate of minimal number of generators}\label{sublude}
Pyber \cite{Py98} conjectures that almost  all finite groups are nilpotent.
It is therefore interesting to get estimates on the minimal number of generators of finite nilpotent groups.

For a prime $p$, a finite group $G$ is called a \emph{$p$-group}
if the orders of all elements of $G$ are divisible by $p$.
This holds if and only if the order of $G$ is a power of $p$;
that is, $|G|=p^\alpha$ for some positive integer $\alpha$.
The number of isomorphism classes $f(n)$ of groups of order $n=p^\alpha$ is given by
\begin{align}\label{higsim}
	f(n) = n^{(2/27 + o(1)) \cdot \alpha^2},
\end{align}
by Higman \cite[Theorem 3.5]{Higman60} and Sims \cite[Proposition 1.1]{Sims65}. 
Because any $p$-group is nilpotent and hence solvable,
the number $f(d,n)$ of groups of order $n=p^{\alpha}$ with a generating set
of at most $d$ elements satisfies
\begin{align}\label{mann}
	f(d, n)  \le n^{(d+1) \cdot \alpha},
\end{align}
by Mann \cite[Theorem 2]{Ma98}.
Combining \eqref{higsim} and \eqref{mann},
we conclude that the proportion of groups of order $n=p^{\alpha}$
with a generating set of at most $d$ elements tends to zero as $\alpha$ tends to infinity.

For a non-trivial finite group $G$, a Sylow $p$-subgroup is a non-trivial subgroup $P$ of $G$ such that $|P|$ is the highest power of $p$ dividing $|G|$.
If $G$ is nilpotent, then $G$ is the product of its $p_i$-Sylow subgroups $P_i$,
\begin{align}
	G = P_1 \times \dots \times P_k,
\end{align}
where the $p_i$ run through the primes dividing $|G|$; see \cite[Theorem 3, Chapter 6]{D-F}.
We conclude

\begin{thm}\label{thmnil}
Let $n_i$ be a sequence of natural numbers such that the maximal exponent of the prime factors of $n_i$ tends to infinity with $i$.
Then the proportion of $d$ generator nilpotent groups of order $n_i$ among all nilpotent groups of order $n_i$ tends to zero with $i$.
\end{thm}

Thus it is asymptotically unlikely that a random nilpotent group of order $n_i$
has a generating set with at most $d$ elements. 
We observe that it is necessary to assume, in \cref{thmnil},
that the maximal exponent of the prime factors of $n_i$ tends to infinity with $i$. 
This follows from a result of Guralnick \cite{Gu89} and Lucchini \cite{Lu89}
that says that a finite group is generated by at most $d+1$ elements
if each of its Sylow subgroups is generated by at most $d$ elements.

Let $M$ be a complete and connected Riemannian manifold of dimension $m$.
Say that two normal Riemannian coverings $p_1\colon M_1\to M$ and $p_2\colon M_2\to M$ are \emph{coarsely equivalent} if their groups of covering transformations are isomorphic.
By \cref{thmnil}, it is unlikely that random coarse equivalence classes of covering transformations of
$n_i$-sheeted normal nilpotent coverings of a surface $S$ have at most $d$ generators as $i$ tends to infinity.
Thus asymtotic $\lambda_1$-stability is unlikely as well as $n_i\to\infty$. 

\section{A basic argument and applications}\label{basic}
Let $\lambda<\lambda_{\ess}$, $\vf\in E_\lambda$,
and $(U_i)_{i\in I}$ be the family of pairwise different nodal domains of $\vf$.
For each $i\in I$, let $(U_{ij})_{j\in J_i}$ be the family of pairwise different nodal domains of $p^*\vf$ over $U_i$.
For each $i\in I$ and $j\in J_i$, let $\vf_{ij}$ be the function on $M'$ which coincides with $p^*\vf/k_{ij}$ on $U_{ij}$, vanishes on all other $U_{ik}$, and is equal to $p^*\vf/|p|$ on the rest of $M'$, where $k_{ij}$ is the degree of the covering $p\colon U_{ij}\to U_i$.
The set $J$ of pairs $ij$ with $i\in I$ and $j\in J_i$ labels the set of all nodal domains $U_{ij}$ of $p^*\vf$, sorted by the nodal domains $U_i$ of $\vf$.
For any function $\psi$ on $M$, we have
\begin{align}\label{integral}
    \int_{M'} \vf_{ij}p^*\psi = \int_M \vf\psi.
\end{align}
Notice the similarity, and difference, between \eqref{padjoint} and \eqref{integral}.
The definition of the $\vf_{ij}$ is adapted to what is needed in the comparison of $N_M$ and $N_{M'}$.
Recall here that $N_M(\lambda)$, $N_{M'}(\lambda)$ and $N_M(\lambda-)$, $N_{M'}(\lambda-)$ denote the number of eigenvalues of $M$ and $M'$ in $[0,\lambda]$ and $[0,\lambda)$, respectively.

\begin{thm}\label{numberg}
There are at least $|J|-|I|$ linearly independent eigenfunctions on $M'$ with eigenvalues in $(0,\lambda)$,
which are perpendicular to $p^*(L^2(M))$.
In particular,
\begin{align*}
    N_{M'}(\lambda-)
    \ge N_M(\lambda-) + \sum_{i\in I}(|J_i| - 1)
    = N_M(\lambda-) + \sum_{i\in I}|J_i| - |I|.
\end{align*}
\end{thm}

Note that $|J_i|\ge1$ for all $i$ so that the summands in the middle,
the \emph{contributions $|J_i| - 1$ of the $U_i$}, are all non-negative.

\begin{proof}[Proof of \cref{thma}]
The contribution of $U$ to the estimate in \cref{numberg} is $j-1$.
Since the contributions of the other nodal domains are non-negative, the claim follows.
\end{proof}

\begin{proof}[Proof of \cref{numberg}]
Let $X$ be the space spanned by the $\vf_{ij}$ and $Y=p^*(\vf^\bot)$.
By \eqref{integral}, $X$ and $Y$ satisfy the assumptions of \cref{respec},
applied to the quadratic form $Q$ associated to the operator $A=\Delta'-\lambda$.
Namely, $Q\le0$ on $X$.
Furthermore, by \eqref{integral}, $X$ and $Y$ are perpendicular in $H=L^2(M')$.
Finally, since $\vf$ is an eigenfunction of $\Delta$, $PY\subseteq Y$,
where $P$ is the spectral projection of $\Delta'$ associated to $[0,\lambda)$.
It remains to clarify the dimensions of $X$ and $X\cap H_0$, where here $H_0=E_\lambda'$.

To determine $\dim X$, suppose that $\sum \alpha_{ij}\vf_{ij}=0$.
Let $i\in I$ and $j\in J_i$.
Then on $U_{ij}$, the $\vf_{ik}$, for $k\ne j$, vanish.
Hence, on $U_{ij}$,
\begin{align*}
    \alpha_{ij}\vf_{ij} = - \sum_{\substack{k\in I\setminus\{i\} \\ l\in J_k}} \alpha_{kl}\vf_{kl}.
\end{align*}
Now on $U_{ij}$, $\vf_{ij}$ is equal to $p^*\vf/k_{ij}$ and each $\vf_{kl}$ on the right to $p^*\vf/|p|$.
Since $p^*\vf\ne0$ on $U_{ij}$, we get
\begin{align*}
    \frac{\alpha_{ij}}{k_{ij}}
    = - \sum_{\substack{k\in I\setminus\{i\} \\ l\in J_k}} \frac{\alpha_{kl}}{|p|}
    =: \alpha_i.
\end{align*}
Hence, on $p^{-1}(U_i)$,
\begin{align*}
    \sum_{j\in J_i} \alpha_{ij}\vf_{ij} = \alpha_i p^*\vf.
\end{align*}
We also get that $\alpha_{ij}=\alpha_ik_{ij}$.
Since $\sum_{j\in J_i}k_{ij}=|p|$ for all $i\in I$, we infer that $\sum_{j\in J_i}\alpha_{ij}=\alpha_i|p|$.
Hence the above displayed equality also holds on the rest of $M'$.
Therefore the $\alpha_i$ satisfy the linear equation
\begin{align*}
    \sum_{i\in I}\alpha_i = 0,
\end{align*}
which has $|I|-1$ independent solutions.
In conclusion,
\begin{align*}
    \dim X = \sum_{i\in I}|J_i| - |I| + 1.
\end{align*}
To determine $\dim X\cap E_\lambda'$, note first that $X\cap E_\lambda'$ contains $p^*\vf$.
Conversely, any linear combination of the $\vf_{ij}$
is a multiple of $p^*\vf$ on any of the nodal domains $U_{kl}$ of $p^*\vf$.
Hence, by the unique continuation property, any smooth function in $X$ is a multiple of $p^*\vf$.
In particular, $X\cap E_\lambda'$ consists of multiples of $p^*\vf$.
Therefore $X\cap E_\lambda'$ has dimension one.
\end{proof}

\begin{cor}\label{simcon}
If $k$ of the nodal domains of $\vf$ are simply connected, then
\begin{align*}
	N_{M'}(\lambda-)\ge N_M(\lambda-) + k(|p|-1).
\end{align*}
\end{cor}

\begin{proof}
The preimage under $p$ of any of the simply connected nodal domains has $|p|$ components.
Hence their contribution to the estimate in \cref{numberg} is $k(|p|-1)$.
Since the contributions of the other nodal domains are non-negative, the claim follows.
\end{proof}

In \cref{incom}, for any nodal domain $U$ of $\vf$,
we let $\Gamma'=p_*\pi_1(M,x')$, for any given $x\in U$ and $x'\in p^{-1}(x)$
(and $\Gamma_U$ as usual).
The assertions are independent of the choice of $x$ and $x'$.
Generalizing \cref{simcon}, we have

\begin{cor}\label{incom}
If $k$ of the nodal domains $U$ of $\vf$ satisfy
\begin{enumerate}
\item
$\Gamma_U\subseteq\Gamma'$, then $N_{M'}(\lambda-) \ge N_M(\lambda-) + k$;
\item
$N_\Gamma(\Gamma_U)\subseteq\Gamma'$, then $N_{M'}(\lambda-) \ge N_M(\lambda-) + k(|p|-1)$.
\end{enumerate}
\end{cor}

\begin{proof}
For any $U$ as in the first assertion, $p^{-1}(U)$ has at least two,
in the second $|p|$ components, by \cref{gugp}.
\end{proof}

\begin{rem}\label{numberc}
Since liftings of eigenfunctions from $M$ to $M'$ are eigenfunctions on $M'$ with the same eigenvalues,
an inequality $N_{M'}(\lambda-)\ge N_M(\lambda-) + C$ implies that $N_{M'}(\kappa-)\ge N_M(\kappa-)+C$
and $N_{M'}(\kappa)\ge N_M(\kappa)+C$, for any $\kappa\ge\lambda$.
\end{rem}

\begin{thm}\label{theoi2}
Suppose that $M$ is complete and connected with $\lambda_1(M)<\lambda_{\ess}(M)$ and carries a $\lambda_1(M)$-eigenfunction $\vf$
such that its nodal set $\mathcal{Z}(\vf)$ has at least $\mu+1\ge2$ components.
Let $U$ be one of the two nodal domains of $\vf$ and $x\in U$.
Then there is a surjective homomorphism $I_x\colon\Gamma\to\Z_2^{\mu}$
such that $\Gamma_U\subseteq\ker I_x$.
If $M$ is orientable, there is a surjective homomorphism $I_x\colon\Gamma\to\Z^{\mu}$
such that $\Gamma_U\subseteq\ker I_x$.
In both cases, if $\Gamma'$ is a finite index subgroup of $\Gamma$ containing $\Gamma_U$,
then the corresponding Riemannian covering $p\colon M'\to M$ is strictly $\lambda_1$-unstable.
More precisely, with $|p|=|\Gamma'\backslash\Gamma|$,
\begin{align*}
	N_{M'}(\lambda_-) \ge N_M(\lambda_-) + \mu(|p| -1).
\end{align*}
\end{thm}

\begin{proof}
Suppose first that $0$ is a regular value of $\vf$.
Then $\mathcal{Z}(\vf)$ is a smooth manifold.
Since $\vf$ is a $\lambda_1(M)$ eigenfunction, it has exactly two nodal domains, $\{\vf<0\}$ and $\{\vf>0\}$.

By assumption, $\mathcal{Z}(\vf)$ has at least $\mu+1$ connected components, $Z_1,\dots,Z_{\mu+1}$.
Let $Z=Z_j$ for some $1\le j\le\mu$, $Z'=Z_{\mu+1}$, and $z\in Z$ and $z'\in Z'$ be points.
Since $\{\vf<0\}$ and $\{\vf>0\}$ are connected,
there exist paths $c_-$ in $\{\vf\le0\}$ and $c_+$ in $\{\vf\ge0\}$ between $z$ and $z'$.
Their union $c$ is a closed loop in $M$ that has intersection number one with $Z$.
In particular, intersection with the different $Z=Z_j$ defines a non-trivial homomorphism $i_Z$
from $\Gamma$ to $\Z^\mu$ if $M$ is oriented and to $\Z_2^\mu$ otherwise.

Let now $U=\{\vf<0\}$ and $x\in U$.
Then $\Gamma_U$ is contained in the kernel $\Gamma'$ of $i_Z$.
Now \cref{gugp} applies and shows the claim.

If $0$ is not a regular value of $\vf$, the above argument still applies in principle,
but we need some preparation to define intersection numbers.
To that end, write
\begin{align*}
	\mathcal{Z}(\vf)
	&= \{z\in\mathcal{Z}(\vf) \mid d\vf(z)\ne0\}
	\cup \{z\in\mathcal{Z}(\vf) \mid d\vf(z)\ne0\} \\ 
	&= \mathcal{Z}(\vf)_{\mathrm{reg}} \cup \mathcal{Z}(\vf)_{\mathrm{sing}},
\end{align*}
the \emph{regular} and \emph{singular part} of $\mathcal{Z}(\vf)$.
Let $Z$ and $Z'$ be as above and set
\begin{align*}
	Z_{\mathrm{reg}} = Z \cap \mathcal{Z}(\vf)_{\mathrm{reg}}, \quad
	Z_{\mathrm{sing}} = Z \cap \mathcal{Z}(\vf)_{\mathrm{sing}}.
\end{align*}
Recall from the (elementary) proof of \cite[Lemma 1.9]{HS89} that any point in $Z_{\mathrm{sing}}$
is contained in an open ball $B$ in $M$ such that $Z_{\mathrm{sing}}\cap B$
is contained in a finite union of embedded submanifolds of dimension $\dim M-2$.
Let now $U=\{\vf<0\}$ and $x\in U$ as above.

Claim 1)
Any loop in $M$ at $x$ is homotopic to a loop which does not meet $Z_{\mathrm{sing}}$
and meets $Z_{\mathrm{reg}}$ transversally in at most finitely many points.

To show 1), let $c$ be a loop in $M$ at $x$.
Then $I=c^{-1}(Z_{\mathrm{sing}})$ is compact.
Hence $I$ can be covered by finitely many consecutive intervals such that the image of each of these intervals is contained in a ball $B$ as above.
Since the codimension of the corresponding embedded submanifolds as above is two,
$c$ can be deformed consecutively to a loop at $x$ which does not meet $Z_{\mathrm{sing}}$.
This shows the first assertion.
The second follows from standard transversality theory, applied to the smooth hypersurface $Z_{\mathrm{reg}}$.

Claim 2)
For any two homotopic loops in $M$ at $x$, which do not meet $Z_{\mathrm{sing}}$
and meet $Z_{\mathrm{reg}}$ transversally in at most finitely many points,
the (oriented respectively non-oriented) intersection numbers with $Z$ coincide.

To show 2), let $c_1$ and $c_2$ be two such homotopic loops.
Since they do not meet $Z_{\mathrm{sing}}$ and intersect $Z_{\mathrm{reg}}$ transversally in at most finitely many points,
there is an $\ve>0$, which is a regular value of $\vf$ such that $\{\vf=\ve\}$ has a component $Z_\ve$
such that the intersections of $c_1$ with $Z$ and $Z_\ve$ are in one-to one correspondence to each other,
and similarly for $c_2$.
Now the intersection numbers of $c_1$ and $c_2$ with $Z_\ve$ are well-defined
and agree with the intersection numbers with $Z_{\mathrm{reg}}$, by what we said.
Since $c_1$ and $c_2$ are homotopic, their intersection numbers with $Z_\ve$ agree,
hence also the ones with $Z_{\mathrm{reg}}$.

Now 1) and 2) show that intersection numbers with the different $Z=Z_j$ respectively $Z_{j,\mathrm{reg}}$
define a homomorphism $\Gamma$ to $\Z^\mu$ in the oriented case and $\Z_2^\mu$ otherwise.
The rest of the proof is as in the regular case.
\end{proof}

\begin{rem}
Under appropriate assumptions, the conclusion of \cref{theoi} holds with $\lambda_k$-unstable in place of $\lambda_1$-unstable.
More precisely, if there is a $\lambda_k(M)$-eigenfunction $\vf$ on $M$ such that there are components of $\mathcal{Z}(\vf)$ together with loops in $M$ which intersect exactly once,
then the above arguments apply and show $\lambda_k$-instability.
This may indicate that $\lambda_k$-instability becomes more likely,
the more nodal domains $\lambda_k(M)$-eigenfunctions have.
\end{rem}

\subsection{Absolute estimate}\label{subdom}
In the above, we estimated $N_{M'}$ against $N_M$.
That is what is behind the definition of the $\vf_{ij}$ in the beginning of the section.
An easier approach leads to an estimate of $N_{M'}$ without comparing it with $N_M$.
The point is as follws:
Let $(D_i)_{i\in I}$ be a family of domains in $M$ and, for each $i\in I$, $(D_{ij})_{j\in J_i}$ be the connected components of $p^{-1}(D_i)$.
Let $\lambda\ge0$ and $\vf{i}$ be a smooth function on $D_i$ with compact support and Rayleigh quotient $\le\lambda$. 
For each $i\in I$ and $j\in J_i$, let $\vf_{ij}$ now be the function which equals $p^*\vf_i$ on $D_{ij}$ and vanishes outside of $D_{ij}$.
Then the $\vf_{ij}$ are pairwise $L^2$-orthogonal and have Rayleigh quotient $\le\lambda$.
Hence
\begin{align}\label{absest}
    N_{M'}(\lambda) \ge \sum_{i\in I}|J_i|.    
\end{align}
Clearly, $N_M(\lambda)\ge|I|$, but that does not lead to an inequality as in \cref{numberg}. 

In the following discussion, we use \eqref{absest} only in the case of one domain in $M$, that is, $|I|=1$; cf. \cref{onedom}.
Set
\begin{align}\label{sigma}
    \sigma(M) = \inf_D\lambda_0(D),
\end{align}
where the infimum is taken over all simply connected domains $D \subset M$.
By monotonicity, $\lambda_0(M)\le\sigma(D)$.
Recall also that $\lambda_0(M')=\lambda_0(M)$.

In general, \eqref{sigma} poses the optimization problem of the existence of a simply connected domain $D$ in $M$ such that $\sigma(M)=\lambda_0(D)$ and of a Dirichlet eigenfunction $\vf{}$ on $D$ for $\lambda_0(D)$. 

\begin{prop}\label{number2}
If $\sigma(M)<\lambda_{\ess}(M)$, then $M'$ has at least $|p|$ eigenvalues in $[0,\sigma(M)]$.
\end{prop}

The assumption $\sigma(M)<\lambda_{\ess}(M)$ is satisfied if $M$ is closed since the essential spectrum of closed Riemannian manifolds is empty.
Recall also that $\lambda_{\ess}(M')=\lambda_{\ess}(M)$ since $p$ is finite.

\begin{proof}[Proof of \cref{number2}]
For any $\sigma(M)<\lambda<\lambda_{\ess}(M)$,
let $D\subseteq M$ be a simply connected domain such that there is a $\vf\in C^\infty_c(M)$ with support in $D$
with Rayleigh quotient $<\lambda$.
There are precisely $|p|$ lifts of $D$ to $M'$, and they are pairwise disjoint.
For any such lift $C$, let $\vf_C\in\ C^\infty_c(M')$ be the function with support in $C$ such that $\vf_C=\vf\circ p$ on $C$.
Then the $\vf_C$s and their gradients are pairwise $L^2$-orthogonal and have the same Rayleigh quotient as $\vf$.
\end{proof}

For $\ell\ge1$, let $\sigma_{\ell}(M)=\inf_D\lambda_0(D)$, where the infimum is taken over all domains $D\subset M$
such that the fundamental group of $D$ is generated by at most $\ell$ elements.
Notice that $\sigma(M)=\sigma_0(M)$.

\begin{prop}\label{numberd}
If the minimal number of generators of $\Gamma'\backslash\Gamma$ is $k$ and $\sigma_\ell(M)<\lambda_{\ess}(M)$ for some $\ell\le k$ ,
then $M'$ has at least $k-\ell+1$ eigenvalues in $[0,\sigma_{\ell}(M)]$.
\end{prop}

\begin{proof}
For any $\sigma_{\ell}(M)<\lambda<\lambda_{\ess}(M)$, there is a domain $D$ in $M$ with $\sigma_{\ell}(M)\le\lambda_0(D)<\lambda$
such that the fundamental group of $D$ is generated by at most $\ell$ elements.
Hence the preimage of $D$ under $p$  has at least $k-\ell+1$ components, by \cref{preu2} and \cref{preu4}.
Therefore $M'$ has at least $k-\ell+1$ eigenvalues in $[0,\lambda)$.
\end{proof}

\begin{rem}
For a Riemannian surface $S$, $\sigma_1(S)$ coincides with the \emph{analytic systol of $S$},
introduced in \cite{BMM16}.
Recall that $S$ has at most $-\chi(S)$ eigenvalues in $[0,\sigma_1(S)]$, by \cite[Theorem 1.5]{BMM17}.
Here we get that the covering surface $S'$ has at least two eigenvalues in $[0,\sigma_1(S)]$,
provided that the minimal number of generators of $\Gamma'\backslash\Gamma$ is at least two.
\end{rem}

\begin{rem}\label{onedom}
As in \cref{numberg}, we can also consider families of pairwise disjoint domains
to get a more general estimate than the one in \cref{numberd}.
Namely, if $(D_i)_{i\in I}$ is a finite family of pairwise disjoint domains in $M$
such that the fundamental group of $D_i$ is generated by at most $\ell_i$ elements
and such that $\lambda=\max\lambda_0(D_i)<\lambda_{\ess}(M)$,
then $M$ has at least $\sum_{i\in I}(k-\ell_i+1)$ eigenvalues in $[0,\lambda)$.
The point is that the different lifts of functions $\vf_i\in C^\infty_c(D_i)$
to the different components of $p^{-1}(D_i)$ are pairwise $L^2$-perpendicular.
\end{rem}

\section{Coverings of surfaces}\label{secsurf}
Let $S$ be a connected surface.
Assume that $S$ is of finite type, that is, $S$ is diffeomorphic to the interior of a compact surface $\bar S$ with boundary.
Equivalently, the Euler characteristic $\chi(S)>-\infty$.
The connected components of $\bar S\setminus S$ consist of circles, called \emph{holes} or \emph{circles at infinity}.

Suppose that $S$ is endowed with a complete Riemannian metric which has a square-integrable eigenfunction $\vf_0$ at the bottom $\lambda_0$ of its spectrum.
Recall that the multiplicity of $\lambda_0$ is one and that $\vf_0$ does not change sign, hence can be chosen to be positive.
If the area $|S|<\infty$, then $\vf_0$ is constant and $\lambda_0=0$.
Eigenfunctions of $S$ for eigenvalues $\lambda>\lambda_0$ are perpendicular to $\vf_0$ and hence change sign.
In particular, they have at least two nodal domains.
The structure of the nodal set of such an eigenfunction $\vf$ was clarified in \cite[Theorem 2.5]{Cheng76}:

\begin{thm}[Cheng]\label{cheng}
The nodal set $\mathcal{Z}(\vf)$ of $\vf$ is a locally finite graph in $S$.
Moreover, $z\in\mathcal{Z}(\vf)$ has valence $2n$
if and only if $\vf$ vanishes to order $n$ at $z$.
The opening angles between the edges at $z$ are equal to $\pi/n$.
Furthermore, $\mathcal{Z}(\vf)$ is a locally finite union of immersed circles and lines.
\end{thm}

We will need the following topological result for the study of nodal domains.

\begin{lem}\label{fintyp}
Let $S$ be a connected surface of finite type and $U\subset S$ an open domain with piecewise smooth boundary.
Assume that the complement $U^c$ of $U$ in $S$ contains only finitely many components which are discs or annuli.
Then $U$ has finite type.
\end{lem}

\begin{proof}
The proof rests on the fact that a surface (orientable or not) has finite type if and only if any family of simple closed curves,
which are not null-homotopic, pairwise disjoint, and pairwise not freely homotopic, is finite;
cf.\,\cite[pp.\;259--260]{Ri63}.
Assuming that $U$ is not of finite type, there is an infinite family $\mathcal F$ of simple closed curves in $U$,
which satisfies these conditions with respect to $U$.

If $c$ is a member of $\mathcal F$ and $c$ is null-homotopic in $S$,
then $c$ bounds an embedded disc $D$ in $S$, $c=\partial D$.
Now $\partial U\cap D$ cannot contain components of $\partial U$ which are line segments
and hence consists of finitely many simple closed curves, which bound discs in $D$ which belong to $U^c$.
There are only finitely many such discs, by assumption.
Hence the union $U'$ of such discs with $U$ is of finite type if and only if $U$ is.
Therefore we can assume from now on that the complement of $U$ does not contain components which are discs.

Let now $c_0$ and $c_1$ be members of $\mathcal F$ which are freely homotopic in $S$.
Then there is an embedded annulus $A$ in $S$ such that $c_0\cup c_1=\partial A$.
Now $\partial U\cap A$ cannot contain components of $\partial U$ which are line segments,
nor can it contain closed curves which are homotopic to zero, by assumption.
Hence it consists of two boundary curves $\hat c_0$ and $\hat c_1$,
such that the parts of $A$ between $c_0$ and $\hat c_0$ respectively $c_1$ and $\hat c_1$ belong to $U$
and the rest, $\hat A$, to $U^c$.
Now $\hat A$ is an annulus.
Hence there are only finitely many such, by assumption.
Since $S$ is of finite type, we arrive at a contradiction to the assumption that $\mathcal F$ is infinite.
\end{proof}

\begin{lem}\label{numnod}
A $\lambda$-eigenfunction $\vf$ on $S$ with $\lambda<\lambda_{\ess}(S)$
has only finitely many nodal domains.
\end{lem}

\begin{proof}
Let $\mathcal N$ be the family of different nodal domains of $\vf$.
For any $U\in\mathcal N$, let $\vf_U$ be the function which coincides with $\vf$ on $U$ and vanishes otherwise.
Then the $\vf_U$ are pairwise $L^2$-perpendicular, are in $H^1(S)$, and have Rayleigh quotient $\le\lambda$.
Since $\lambda<\lambda_{\ess}(S)$,
the variational characterization of eigenvalues implies that $\mathcal N$ is finite.
\end{proof}

\begin{cor}\label{geonod}
For $\lambda<\lambda_{\ess}(S)$,
the nodal domains of $\lambda$-eigenfunctions on $S$ are geometrically finite.
\end{cor}

\begin{proof}
By \cref{cheng}, any nodal domains $U$ of $\vf$ is a domain in $S$ with piecewise smooth boundary.
(If $\partial U$ has a self-intersection at a critical point $x$ of $\vf$,
push $U$ a bit inside, away from $x$, to get $\partial U$ embedded.)
Now \cref{fintyp} and \cref{numnod} imply the assertion.
\end{proof}

\subsection{On the topology of nodal domains}\label{subtop}
Assume from now on that $\chi(S)\le0$ and let $\vf$ be an eigenfunction of $S$ perpendicular to $\vf_0$.

\begin{lem}\label{fintyp2}
If none of the nodal domains of $\vf$ is a disc, then each nodal domain $U$ of $\vf$ is a domain of finite type; in particular $\chi(U)>-\infty$.
Moreover $\chi(U)<0$ unless $U$ is a disc or an annulus or a M\"obius band.
\end{lem}

\begin{proof}
Let $U$ be a nodal domain of $\vf$.
If a component $D$ of $U^c$ is a closed disc, then a component of the open set $D\setminus\mathcal{Z}(\vf)$ is a disc,
a contradiction to the assumption.
Now \cref{fintyp} implies the assertion.
\end{proof}

\begin{lem}\label{nodalchi}
If $\vf$ is an eigenfunction of $S$ perpendicular to $\vf_0$, then $\vf$ has a nodal domain $U$ such that $\chi(S)/2\le\chi(U)\le 1$.
More generally, if $\ell\ge2$ denotes the number of nodal domains of $\vf$,
then $\vf$ has a nodal domain $U$ such that
\begin{align*}
	\chi(S)/\ell \le \chi(U) \le 1.
\end{align*}
\end{lem}

\begin{proof}
Since $\chi(S)\le0$, we can assume that no nodal domain of $\vf$ is a disc.
Then all nodal domains of $\vf$ are domains of finite type, in particular with finite Euler characteristic.
By \cite[Theorem 2.5]{Cheng76}, the nodal set $Z$ of $\vf$ is a graph with vertices of (even) order at least two. 
Therefore $\chi(Z)\le0$ and hence
\begin{align*}
	\sum \chi(U) \ge \chi(Z) + \sum \chi(U) = \chi(S),
\end{align*}
where the sum is over all nodal domains $U$ of $\vf$.
Hence there is at least one nodal domain $U$ of $\vf$ such that $\chi(S)/\ell\le\chi(U)\le 1$.
\end{proof}

Let $\Gamma$ be the fundamental group of $S$.
If $S$ is closed, $S=\bar S$, the minimal number of generators of $\Gamma$ is
\begin{align}\label{k2}
    \nu = \nu(\Gamma) = \nu(S) = 2-\chi(S).
\end{align}
If $S$ is non-compact, then $\Gamma$ is a free group and the minimal number of generators of $\Gamma$ is
\begin{align}\label{k1}
    \nu = \nu(\Gamma) = \nu(S) = 1-\chi(S).
\end{align}

\begin{lem}\label{nodalfun}
Let $\vf$ be an eigenfunction of $S$ with $\ell\ge2$ nodal domains.
\begin{enumerate}
\item
If $\ell\ge\nu(S)$, at least one of the nodal domains is a disc or an annulus or a M\"obius band.
\item
If $\ell<\nu(S)$, then $\vf$ has a nodal domain with minimal number of generators of its fundamental group at most $1-\chi(S)/\ell$.
\end{enumerate}
\end{lem}

\begin{proof}
If all nodal domains of $\vf$ have negative Euler characteristic, then $\ell\le-\chi(S)$.
The first claim now follows from \eqref{k2} and \eqref{k1}.
As for the second, we may assume that all nodal domains of $\vf$ have negative Euler characteristic.
Hence, by \cref{nodalchi}, $\vf$ has a nodal domain $U$ with $-\chi(S)/\ell\ge-\chi(U)$.
But then $\nu(U) = 1-\chi(U) \le 1-\chi(S)/\ell$ by \eqref{k1}.
\end{proof}

\begin{cor}\label{nodalfunc}
There is a nodal domain $U$ of $\vf$ such that the fundamental group of $U$ admits a system with at most $\nu/2$ generators if $S$ is closed and $(\nu+1)/2$ generators otherwise, where $\nu=\nu(S)$.
\end{cor}

Consider now a finite Riemannian covering $p\colon S'\to S$ of complete and connected Riemannian surfaces of finite type.
Write $S=\Gamma\backslash\tilde S$ and  $S'=\Gamma'\backslash\tilde S$,
where the fundamental groups $\Gamma\supseteq\Gamma'$ of $S$ and $S'$
are viewed as groups of covering transformations of the universal covering surface $\tilde S$ of $S$ and $S'$.

\begin{prop}\label{nonana}
If $p$ is $\lambda_1$-stable, where $\lambda_1(S)<\lambda_{\ess}$,
then the minimal number of generators of $\Gamma'\backslash\Gamma$ is at most $\nu/2$ if $S$ is closed and $(\nu+1)/2$ otherwise, where $\nu=\nu(S)$. 
\end{prop}

\begin{proof}
Let $\vf$ be a $\lambda_1(S)$-eigenfunction.
By \cref{nodalfunc}, $\vf$ has a nodal domain $U$ such that the fundamental group of $U$ admits a system with at most $\nu/2$ respectively $(\nu+1)/2$ generators.
If the minimal number of generators for $\Gamma'\backslash\Gamma$ is strictly bigger than $\nu/2$ respectively $(\nu+1)/2$,
then the right action of $\pi_1(U)$ on $\Gamma'\backslash\Gamma$ cannot be transitive.
By \cref{preu2} together with \cref{numberg}, we get a contradiction to $\lambda_1$-stability.
\end{proof}

\subsection{Intersection numbers and coverings}\label{subinter}
Let $U$ be a nodal domain of an eigenfunction $\vf$ of $S$ perpendicular to $\vf_0$.
Then the complement $U^c$ of $U$ is a surface of finite type,
and we let $V$ be a component of $U^c$.
Then $V$ has $k\ge1$ boundary circles in $S$, that we call \emph{doors} $d_1,\dots,d_k$, through which it is connected to $U$
and $\ell\ge0$ boundary circles in $\bar S\setminus S$, that we call \emph{exits $e_1,\dots,e_\ell$ at infinity}.
We draw the first kind in green, meaning that we may enter $U$ through them, and the latter kind in red, indicating that we exit $S$ through them eventually.

We now view $V$ as a regular plane polyhedron $P$ in the orientable case, respectively $Q$ in the non-orientable case, with the colored holes in its interior
and with the standard identifications of its edges, indicated by the labellings
\begin{alignat*}{2}
	&P_0\colon aa^{-1}, &\quad\text{where $g=0$,} \\
	&P_g\colon a_1b_1a_1^{-1}b_1^{-1}\cdots a_gb_ga_g^{-1}b_g^{-1}, &\quad\text{where $g\ge1$,}  \\
	&Q_g\colon a_1a_1\cdots a_ga_g, &\quad\text{where $g\ge1$.} 
\end{alignat*}
In terms of (orientable respectively non-orientable) \emph{genus} $g$ and numbers $k$ and $\ell$ of holes,
the negative of the Euler characteristic of $V$ is
\begin{align}\label{chiv}
	-\chi(V) = 2g + k  + \ell - 2 \quad\text{and}\quad -\chi(V) = g + k + l -2
\end{align}
in the orientable ($P_g$) and non-orientable ($Q_g$) case, respectively.

 If $\ell\ge2$, we draw disjoint segments $h_2,\dots,h_\ell$ from $e_1$ to the other red circles $e_2,\dots,e_\ell$. 
We get a homomorphism
\begin{align}
	I = I_V\colon H_1(S) \to \Z^{\mu}
	\quad\text{respectively}\quad
	I = I_V\colon H_1(S) \to \Z_2^{\mu}, 
\end{align}
the \emph{intersection homomorphism}, where
\begin{align}\label{ihom}
\begin{alignedat}{2}
	\mu = \mu(V) &= 2g + k -1 + (\ell-1)\vee0 & &\text{for $P_g$,} \\
	\mu = \mu(V) &= g + k -1 + (\ell-1)\vee0 & \quad&\text{for $Q_g$,}
\end{alignedat}
\end{align}
by taking,
\begin{enumerate}
\item\label{Io}
in the oriented case, oriented intersection numbers with the circles in $V$ coming from the edges of $P$ labeled $a_1,b_1,\dots,a_g,b_g$, the boundary circles $d_2\dots,d_k$, and the segments $f_2\dots,f_\ell$.
\item\label{Ino}
in the non-orientable case, intersection numbers modulo two with the circles in $V$ coming from the edges of $Q$ labeled $a_1,\dots,a_g$, the boundary circles $d_2\dots,d_k$, and the segments $f_2\dots,f_\ell$.
\end{enumerate}

From \eqref{chiv} and \eqref{ihom}, we conclude

\begin{lem}\label{chim}
Irrespective of whether $V$ is orientable or not, we have
\begin{align*}
	\mu(V) = - \chi(V) + \begin{cases} 1 &\text{if $\ell=0$} \\ 0 &\text{if $\ell\ge1$}\end{cases}
	\ge - \chi(V).
\end{align*}
\end{lem}

\begin{rem}\label{m01}
We have $\mu=\mu(V)\ge1$ unless $(g,k,\ell)$ equals $(0,1,0)$ or $(0,1,1)$,
and then $V$ is a disc with one green or an annulus with one green and one red boundary circle, respectively.
There are also only a few cases with $\mu=1$.
Using the above $(g, k, l)$ notation, if $V$ is orientable, then it is one of the following:
\begin{enumerate}
\item[(0, 2, 0)]\label{020}
an annulus, boundary circles green; $\chi(V)=0$,
\item[(0, 2, 1)]\label{021}
a pair of pants, two boundary circles green, one red,
\item[(0, 1, 2)]\label{012}
a pair of pants, one boundary circle green, two red.
\end{enumerate}
If $V$ is non-orientable, it is one of the following:
\begin{enumerate}[resume]
\item[(1, 1, 0)]\label{110}
a \emph{M\"obius band}, boundary circle green; $\chi(V)=0$.
\item[(1, 1, 1)]\label{111}
a projective plane with two holes, boundary circles green and red.
\end{enumerate}
\end{rem}

\begin{lem}\label{gukiv}
Under the Hurewicz homomorphism (see \cref{abco}) from $\Gamma$ to $H_1(S)$ respectively $H_1(S;\Z_2)$,
$\Gamma_U$ is contained in the kernel $\ker I_V$.
\end{lem}

\begin{proof}
By definition, homotopy classes of loops in $\Gamma_U$ have representatives which are contained in $U$
and hence have empty intersection with curves in $V\subseteq U^c$. 
\end{proof}

\begin{thm}\label{thmbk}
Let $p\colon S'\to S$ be a finite Riemannian covering of complete and connected Riemannian surfaces
and $\vf$ be a $\lambda$-eigenfunction of $M$, where $\lambda_0(S)<\lambda<\lambda_{\ess}(S)$.
Let $k$ be the number of nodal domains $U$ of $\vf$ such that, for some component $V$ of $U^c$, $\ker I_V\subseteq\Gamma'$,
where $\Gamma'$ is the image of $\pi_1(S',x')$ in $\pi_1(S, x)$ under $p_\#$ for some $x'\in p^{-1}(x)$ (see \cref{secco}).
Then
\begin{align*}
	N_{S'}(\lambda-) \ge N_S(\lambda-) + k(|p|-1).
\end{align*}
\end{thm}

\begin{proof}
By \cref{gukiv} and since $\ker I_V$ is a normal subgroup of $\Gamma$,
$N_{\Gamma}(\Gamma_U)$ is contained in $\ker I_V$,
for any nodal domain $U$ and any component $V$ of its complement as in the assertion.
Hence $p^{-1}(U)$ has $|p|$ components for any such $U$, by \cref{gugp}.
Now \cref{numberg} implies the assertion.
\end{proof}

\begin{thm}\label{thmb2}
Assume that $S$ is of finite type with $\chi(S)<0$, and let $\vf$ be a $\lambda$-eigenfunction,
where $\lambda_0(S)<\lambda<\lambda_{\ess}(S)$.
Then $\vf$ has $\nu\ge2$ nodal domains and at least one, $U$, such that $\chi(S)/\nu\le\chi(U)\le1$.
For any such $U$, $\Gamma$ admits a surjective homomorphism $I$ to $\Z^{\mu}_2$,
respectively $\Z^{\mu}$ if $S$ is orientable,
where $\mu\ge-(\nu-1)\chi(S)/\nu$, such that $\Gamma_U\subseteq\ker I$.
In particular, if $\ker I\subseteq\Gamma'\subseteq\Gamma$ is a finite index subgroup,
then the corresponding Riemannian covering $p\colon S'\to S$ is strictly $\lambda$-unstable.
More precisely,
\begin{align*}
	N_{S'}(\lambda-) \ge N_S(\lambda-) + |p| -1,
\end{align*}
where $|p|=|\Gamma'\backslash\Gamma|$.
\end{thm}

\begin{proof}
By \cref{nodalfunc}, $\vf$ has a nodal domain $U$ such that $\chi(U)\ge\chi(S)/\nu$.
Then, by \cref{cheng},
\begin{align*}
	\chi(U^c) = \chi(S) - \chi(U) \le (\nu-1)\chi(S)/\nu < 0.
\end{align*}
Therefore the components $V_j$ of $U^c$ with $\chi(V_j)<0$ satisfy
\begin{align}
	\sum_j\chi(V_j) \le (\nu-1)\chi(S)/\nu.
\end{align}
For each $i$, we have
\begin{align*}
	\mu(V_j) \ge - \chi(V_j),
\end{align*}
by \cref{chim}.
Therefore, if $\mu$ equals the sum of the $\mu(V_j)$, then the sum
\begin{align*}
	I = \oplus_j I_{V_j} \colon \oplus_j H_1(V_j) \to \oplus_j \im I_{V_j}
\end{align*}
is a homomorphism to $\Z^{\mu}$ respectively $\Z_2^{\mu}$ as asserted,
except that we compose it with the corresponding Hurewicz homomorphism
to have it defined on $\Gamma$.
\end{proof}

\subsection{Estimating the number of unstable coverings}\label{subcou}
Fix a base point $x\in S$,
and consider finite Riemannian coverings $p\colon(S',x')\to(S,x)$ of pointed complete and connected Riemannian surfaces.
Note that the isomorphism classes of such pointed coverings with $n$ sheets correspond one-to-one with index $n$ subgroups of $\Gamma=\pi_1(S,x)$.
Denote by $a(n)$ the number of isomorphism classes of all such $n$-sheeted coverings and by $u(n)$ the number of isomorphism classes of $\lambda_1$-unstable ones among them.

\begin{cor}\label{count}
If $S$ is of finite type with $\chi(S)<0$, then \[u(2n)\ge \frac{a(n)}{2n-1}.\]
\end{cor}

\begin{proof}
We want to estimate $a(n)$ against $u(2n)$.
Now for any $n$-sheeted pointed covering $p\colon(S',x')\to(S,x)$ as above, \cref{thmb2} applies to $S'$ in place of $S$
and shows, that there is a $\lambda_1$-unstable twofold Riemannian covering  \[p'\colon(S'',x'')\to(S',x').\]
Then the composition $q=p\circ p'$ is $2n$-sheeted,
and $q$ is $\lambda_1$-unstable, whether $p$ is $\lambda_1$-unstable or not.
Now the number of isomorphism classes of such $p$ versus the given $q$
is estimated by an upper bound on the number of index $n$ subgroups $\Gamma'$ of $\Gamma$ containing the given $\Gamma''$ of index $2n$.
Since the index of $\Gamma''$ in $\Gamma'$ is two,
$\Gamma'$ is generated by $\Gamma''$ together with an element $g\in\Gamma\setminus\Gamma''$.
Furthermore, since the index of $\Gamma''$ in $\Gamma$ is $2n$, there are at most $2n-1$ such $g$ modulo $\Gamma''$
(such that the subgroup $\Gamma'$ generated by $g$ and $\Gamma''$ has index $n$ in $\Gamma$).
\end{proof}

\begin{rem}\label{count2}
Let $\Gamma=\pi_1(S,x)$ and $\tilde S$ be the universal covering surface of $S$, endowed with the lifted metric.
Let $\rho$ be a homomorphism from $\Gamma$ to the symmetric group $S_n$.
Associated to $\rho$, Magee et al.\;consider the orbit space $S'$ of the product action of $\Gamma$ on $\tilde S\times\{1,\dots,n\}$,
see \cite[paragraph following (1.1)]{MageeNaudPuder22}.
The natural projection $p\colon S'\to S$ is an $n$-sheeted covering,
and the construction also yields a labeling of $p^{-1}(x)$ by $\{1,\dots,n\}$, a \emph{labeled covering}.
Notice that $S'$ is connected if and only if the action of $\rho(\Gamma)$ on $\{1,\dots,n\}$ is transitive.

Clearly, isomorphism classes of labeled coverings of $(S,x)$ are characterized by the representations $\rho$.
Disregarding the labeling, but fixing a base point $x'\in p^{-1}(x)$ corresponds to dividing the labeling by $S_{n-1}$. 
Thus the number of all isomorphism classes of $n$-sheeted connected and labeled coverings equals $(n-1)!a(n)$ with $a(n)$ as above. 
Since $\lambda_1$-stability is independent of the choice of base points and labelings,
we conclude that the analog of the inequality of \cref{count} holds for isomorphism classes of connected labeled coverings as well. 

Magee et al.\;point out that, as $n\to\infty$,
the number of non-connected $n$-sheeted labeled coverings vanishes asymptotically in proportion to the connected ones.
\end{rem}

\section{Towers of coverings}\label{subman}
Let $\cdots \to M_k \to M_{k-1} \to \dots \to M_1 \to M_0 = M$ be a tower of Riemannian coverings of complete and connected Riemannian manifolds.
Assume that the degrees $|p_k|$ of the coverings $p_k\colon M_k\to M$ tend to infinity as $k$ tends to infinity.

We say that a sequence of points $x_k\in M_k$ is a \emph{tower of points} above $x=x_0\in M_0$ if they lie above each other.
Given $x\in M$, we choose a tower $(x_k)$ of points above $x$
and get the fundamental groups $\Gamma=\pi_1(M,x)$ and $\Gamma_k=\pi_1(M_k,x_k)$.
By choosing a point $\tilde x$above the $x_k$ in the universal covering space $\tilde M$,
we identify $\Gamma$ and the $\Gamma_k$ with the corresponding groups of covering transformations of $\tilde M$ such that
\begin{align*}
    \cdots \subseteq \Gamma_k \subseteq \Gamma_{k-1} \subseteq \dots \subseteq \Gamma_1 \subseteq \Gamma_0 = \Gamma
\end{align*}
and $M_k=\Gamma_k\backslash\tilde M$,  $M=\Gamma\backslash\tilde M$.

\begin{lem}\label{diam}
For any $x_0\in M$, $\diam p_k^{-1}(x_0)\to\infty$ as $k\to\infty$.
\end{lem}

\begin{proof}
Let $r>0$ be the injectivity radius of $M$ at $x_0$.
Then the geodesic ball $B(x_0,r)$ is evenly covered by any Riemannian covering of $M$.

Let $\tilde p\colon\tilde M\to M$ be the universal covering.
Endowed with the lifted Riemannian metric, $\tilde M$ is complete and the coverings $\tilde M\to M_k$ are Riemannian.

Let $k\ge1$ and $y_1,\dots,y_l$ be the points in $p_k^{-1}(x_0)$, where $l=|p_k|$.
Then $\diam p_k^{-1}(x_0)$ is realized by the minimal $\diam\{z_1,\dots,z_l\}$,
where the minimum is taken over all choices of lifts $z_i$ of the $y_i$ to $\tilde M$.
The minimum is attained since $\tilde p^{-1}(x_0)$ is a discrete subset of $\tilde M$.
Now the assertion follows by the same reason and since  $|p_k|\to\infty$ as $k\to\infty$.
\end{proof}

For a tower of Riemannian coverings as above, we say that a covering $\hat p\colon\hat M\to M$ is a \emph{roof} of the tower
if $\hat p$ is a Riemannian covering of complete and connected Riemannian manifolds which factors through the $p_k$.
The universal covering of $M$ is the highest roof.
The lowest roof is given by $\cap\Gamma_k$. 

The following result extends \cref{tower}.

\begin{thm}\label{roof}
For a tower of finite Riemannian coverings together with a roof $\hat M$ as above,
suppose that the covering $\hat p\colon\hat M\to M$ is normal and that $\lambda_0(\hat M)<\lambda_{\ess}(M)$.
Then, for any $\lambda_0(\hat M)<\lambda<\lambda_{\ess}(M)$ and $l\ge1$,
$M_k$ has at least $l$ eigenvalues below $\lambda$, for all sufficiently large $k$; in short,
\[\limsup_{k\to\infty}\lambda_l(M_k)\le\lambda_0(\hat M).\]
\end{thm}

\begin{proof}
For $\lambda$ as in the assertion,
there is a smooth function $\vf$ on $\hat M$ with support contained in some geodesic ball $B(x,r)\subset\hat M$
and with Rayleigh quotient $\Ray\vf<\lambda$.

Let $x_0\in M$ be the point under $x$.
Since the degrees of the coverings $p_k\colon M_k\to M$ tend towards $\infty$,
the same is true for the diameter of the fibers $F_k$ of $p_k$ over $x_0$, by \cref{diam}.
In particular, there is a first $k_2\ge1$ such that $F_{k_2}$ contains two points $y_{k_2,1},y_{k_2,2}$ of distance $d_2\ge2r$.
Then any fiber $F_k$ with $k\ge k_2$ contains two points $y_{k,1},y_{k,2}$ of the same distance $d_2$.
There is then a first $k_3\ge k_2$ such that $F_{k_3}$ contains a further point $y_{k_3,3}$ of distance $\ge2r$ to $y_{k_3,1}$ and $y_{k_3,2}$.
Then any fiber $F_k$ with $k\ge k_3$ contains a further point $y_{k,3}$ with distances
\begin{align*}
    d_3=d(y_{k,1},y_{k,3}) = d(y_{k_3,1},y_{k_3,3})
    \quad\text{and}\quad d(y_{k,2},y_{k,3}\ge2r.
\end{align*}
Proceeding in this way, given any $l\ge1$, there is $k_l\ge1$ such that $F_k$, for $k\ge k_l$,
contains $l$ points of pairwise distance $\ge2r$.
The pushdown, in the sense of \cite[Section 4]{BMP18},
of translates of $\vf$ to preimages of these points in $\hat M$
then yields functions with compact support in the balls of radius $r$ about these points and with Rayleigh quotients $<\lambda$.
Since the supports of these functions in $M_k$ are disjoint, they and their gradients are pairwise $L^2$-orthogonal.
Hence $\lambda_l(M_k)<\lambda$, by the variational characterization of eigenvalues.
\end{proof}

\begin{cor}\label{toweram}
Under the assumptions of \cref{tower}, if the group of covering transformations of $\hat p$ is amenable,
then $\lim_{k\to\infty}\lambda_l(M_k)=\lambda_0(M)$, for any $l\ge1$.
\end{cor}

\begin{proof}
Since the group of covering transformations of $\hat p$ is amenable,
we have $\lambda_0(\hat M)=\lambda_0(M)$, by \cite[Theorem 1.2]{BMP18}.
The asserted equality follows now from $\lambda_0(M_k)\ge\lambda_0(M)$,
which is a consequence of the characterization of $\lambda_0$ as the supremum of the positive spectrum.
\end{proof}

\appendix

\section{A remark about spectral theory}
Let $A$ be an unbounded operator with dense domain $\mathcal{D}$ in a Hilbert space $H$.
Consider the polar decomposition $A=U|A|$ of $A$ and the associated orthogonal decomposition
\begin{align}\label{deco}
    H = H_- \oplus H_0 \oplus H_+
\end{align}
with $Ux=\pm x$ for $x\in H_\pm$ and $H_0=\ker A=\ker U$ as in \cite[Section VI.7]{Kato80}.
Since $A$ and $|A|$ commute with $U$, the decomposition of $H$ is invariant under $A$ and $|A|$.
In particular,
\begin{align*}
    \mathcal{D} = (\mathcal{D}\cap H_-) \oplus H_0 \oplus (\mathcal{D}\oplus H_+)
\end{align*}
and, similarly,
\begin{align*}
    \mathcal{D}_Q = (\mathcal{D}_Q\cap H_-) \oplus H_0 \oplus (\mathcal{D}_Q\oplus H_+)
\end{align*}
for the domain $\mathcal{D}_Q\supseteq\mathcal{D}$ of the quadratic form $Q=Q(x,y)=\la Ax,y\ra$ in $H$ associated to $A$.
Since $\ker A=\ker|A|=H_0$, we have $Q<0$ on $\mathcal{D}_Q\cap H_-$ and $Q>0$ on $\mathcal{D}_Q\cap H_+$.

The following is a refined version of the usual variational characterization of eigenvalues of $A$ in the case where the spectrum of $A$ is discrete.

\begin{lem}\label{respec}
Let $X\subseteq\mathcal{D}_Q$ and $Y\subseteq H$ be subspaces such that $Q\le0$ on $X$, $X\perp Y$, and $PY\subseteq Y$,
where $P$ denotes the orthogonal projection of $H$ onto $H_-$.
Then
\begin{align*}
    X\cap\ker P=X\cap H_0 \quad\text{and}\quad PX\perp Y.
\end{align*}
In particular, if $\dim X<\infty$, then
\begin{align*}
    \dim (H_-\ominus Y) \ge \dim PX = \dim X - \dim(X\cap H_0).
\end{align*}
\end{lem}

\begin{proof}
Write $x\in X$ as  $x=x_-+x_0+x_+$ according to \eqref{deco}, where $Px=x_-$.
Now we have
\begin{align*}
    Q(x,x) = Q(x_-,x_-) + Q(x_+,x_+) \le 0
\end{align*}
and hence
\begin{align*}
    Q(x_+,x_+) \le -Q(x_-,x_-).
\end{align*}
Since $Q>0$ on $H_+$, $x_-=0$ implies that $x_+=0$ so that then $x=x_0\in H_0$.
This shows the first assertion.
As for the second, we have
\begin{align*}
    \la PX,Y\ra = \la X,PY\ra = 0
\end{align*}
since $P$ is orthogonal, $PY\subseteq Y$, and $X\perp Y$.
\end{proof}




\begin{thebibliography}{B}

\bibitem{BMM16}
W.\,Ballmann, H.\,Matthiesen, and S.\,Mondal,
Small eigenvalues of closed surfaces.
\emph{Journal of Differential Geometry} \textbf{103} (2016), no. 1, 1--13.

\bibitem{BMM17}
W.\,Ballmann, H.\,Matthiesen, and S.\,Mondal,
Small eigenvalues of surfaces of finite type.
\emph{Compos. Math.}, \textbf{153} (2017), no. 8, 1747--1768.

\bibitem{BMM18}
W.\,Ballmann, H.\,Matthiesen, and S.\,Mondal,
Small eigenvalues of surfaces: old and new.
\emph{ICCM Not.} \textbf{6} (2018), no. 2, 9--24.

\bibitem{BMP18}
W.\,Ballmann, H.\,Matthiesen, and P.\,Polymerakis,
On the bottom of spectra under coverings.
\emph{Math. Z.} \textbf{288} (2018), no. 3--4, 1029--1036.

\bibitem{BP21}
W.\,Ballmann and P.\,Polymerakis,
Bottom of spectra and coverings of orbifolds.
\emph{Internat. J. Math.} \textbf{32} (2021), no. 12, Paper No. 2140018, 33 pp.



\bibitem{Brooks86}
R.\,Brooks,
The spectral geometry of a tower of coverings.
\emph{J. Differential Geom.} \textbf{23} (1986), no. 1, 97--107. 

\bibitem{Cheng76}
S.\,Y.\,Cheng,
Eigenfunctions and nodal sets. \emph{Comment. Math. Helv.} \textbf{51} (1976), no. 1, 43--55.

\bibitem{D-F} 
D.\,S.\,Dummit and R.\,M.\,Foote,
\emph{Abstract Algebra}. Third Edition.
John Wiley \& Sons, Inc., Hoboken, NJ, 2004, xii+932.

\bibitem{EP15}
A.\,Enciso and D.\,Peralta-Salas,
Eigenfunctions with prescribed nodal sets.
\emph{J. Diff. Geom.} \textbf{101} (2015), no. 2, 197-211.

\bibitem{Gu89}
R.\,M.\,Guralnick,
On the number of generators of a finite group.
\emph{Arch. Math.} \textbf{53} (1989), 521-523.

\bibitem{HS89}
R.\,Hardt and L.\,Simon,
Nodal sets for solutions of elliptic equations.
\emph{J. Differential Geom.} \textbf{30} (1989), no. 2, 505--522.

\bibitem{Higman60}
H.\,Higman, 
Enumerating p-groups, I.
\emph{Proc. Lond. Math. Soc.} \textbf{10} (1960), 24--30.

\bibitem{HildeMagee21}
W.\,Hilde and M.\,Magee,
Near optimal spectral gaps for hyperbolic surfaces.
\emph{Ann. of Math.} (2) \textbf{198} (September 2023), no. 2, 791--824. 

\bibitem{Kato80}
T.\,Kato,
{\em Perturbation theory for linear operators.} Reprint of the 1980 edition.
Classics in Mathematics. Springer-Verlag, Berlin, 1995. xxii+619 pp.

\bibitem{Ko06}
R.\,Komendarczyk,
On the contact geometry of nodal sets.
\emph{Trans. Amer. Math. Soc.} \textbf{358} (2006), no. 6, 2399--2413.



\bibitem{Lu89}
A.\,Lucchini, A bound for the number of generators of a finite group.
\emph{Arch. Math.} \textbf{53} (1989), 313-317.

\bibitem{MageeNaud20}
M.\,Magee and F.\,Naud,
Explicit spectral gaps for random covers of Riemann surfaces.
\emph{Publ. Math. Inst. Hautes \'Etudes Sci.} \textbf{132} (2020), 137--179.

\bibitem{MageeNaudPuder22}
M.\,Magee, F.\,Naud, and D.\,Puder,
A random cover of a compact hyperbolic surface has spectral gap {$\frac{3}{16}-\varepsilon$}.
\emph{Geom. Funct. Anal.} \textbf{32} (2022), no. 3, 595--661.

\bibitem{Ma98}
A.\,Mann,
Enumerating finite groups and their defining relations.
\emph{J. Group Theory} \textbf{1} (1998), no.1, 59--64.

\bibitem{Mondal15}
S.\,Mondal,
On largeness and multiplicity of the first eigenvalue of finite area hyperbolic surfaces.
\emph{Math. Z.} \textbf{281} (2015), no. 1-2, 333–348. 


\bibitem{Py98}
L.\,Pyber,
Group enumeration and where it leads us.
\emph{European {C}ongress of {M}athematics, {V}ol. {II} {B}udapest, 1996},
Progr. Math. 169, 187--199, Birkh\"{a}user, Basel 1998.

\bibitem{Ra74}
B. Randol,
Small eigenvalues of the Laplace operator on compact Riemann surfaces.
\emph{Bull. Amer. Math. Soc.} \textbf{80} (1974), 996--1000.

\bibitem{Ri63}
I.\,Richards,
 On the classification of noncompact surfaces.
\emph{Trans. Amer. Math. Soc.} \textbf{106} (1963), 259--269.

\bibitem{Sims65}
C.\,C.\,Sims, 
Enumerating p-groups, \emph{Proc. Lond. Math. Soc.} 15 (1965), 151-166.

\bibitem{Su87}
D.\,Sullivan,
Related aspects of positivity in {R}iemannian geometry.
\emph{J. Differential Geom.},
\textbf{25} (1987), no. 3, 327--351.

\bibitem{Sunada85}
T.\,Sunada,
Riemannian coverings and isospectral manifolds.
\emph{Ann. of Math.} (2) \textbf{121} (1985), no. 1, 169--186. 

\end{thebibliography}
\end{document}